%% file: l4dc2021-sample.tex
\documentclass[12pt]{l4dc2021} 

\title[Submodular Maximizatoin with Distributed Constraints]{Optimal Algorithms for Submodular Maximization with Distributed Constraints}
\usepackage{times}

\author{%
 \Name{Alexander Robey} \Email{arobey1@seas.upenn.edu}
 \AND
 \Name{Arman Adibi} \Email{aadibi@seas.upenn.edu}
 \AND
 \Name{Brent Schlotfeldt} \Email{brentsc@seas.upenn.edu} 
 \AND
 \Name{Hamed Hassani} \Email{hassani@seas.upenn.edu}
 \AND
 \Name{George J.\ Pappas} \Email{pappasg@seas.upenn.edu}\\
 \addr University of Pennsylvania, Philadelphia, PA 19103, USA
}

\makeatletter
 \let\Ginclude@graphics\@org@Ginclude@graphics 
\makeatother

\input{commands}

\begin{document}

\maketitle

\input{contents/abstract}


\input{contents/intro}
\input{contents/related-work}
\input{contents/prelims}
\input{contents/prob-statement}
\input{contents/algorithm}
\input{contents/convergence}
\input{contents/simulation}
\input{contents/conclusion}

\newpage

\bibliography{bib}

\newpage

\begin{appendix}
    \input{appendix/assumptions}
    \input{appendix/prelim-lemmas}

    \input{appendix/main-thm-proof}
\end{appendix}

\end{document}

%% file: commands.tex
\usepackage{graphicx}
\usepackage{xcolor}
\usepackage{enumitem}
\usepackage{algorithm}
\usepackage{algorithmic}

\newcommand{\Z}{\mathbb{Z}}
\newcommand{\R}{\mathbb{R}}
\newcommand{\W}{\mathbf{W}}
\newcommand{\1}{\mathbf{1}}
\newcommand{\0}{\mathbf{0}}
\newcommand{\Y}{\mathcal{Y}}
\newcommand{\X}{\mathcal{X}}
\newcommand{\A}{\mathcal{A}}
\newcommand{\B}{\mathcal{B}}
\newcommand{\I}{\mathcal{I}}
\newcommand{\G}{\mathcal{G}}
\newcommand{\E}{\mathcal{E}}
\newcommand{\setS}{\mathcal{S}}
\newcommand{\setP}{\mathcal{P}}
\newcommand{\setC}{\mathcal{C}}
\newcommand{\setJ}{\mathcal{J}}
\newcommand{\bigO}{\mathcal{O}}
\newcommand{\agents}{\mathcal{N}}
\newcommand{\yregcon}{\yreg_{\text{con}}}
\newcommand{\ybarcon}{\ybar_{\text{con}}}
\newcommand{\vregcon}{\vreg_{\text{con}}}
\newcommand{\identity}{\mathbf{I}}
\newcommand{\indicator}{1}

\newtheorem{problem}[theorem]{Problem}
\newtheorem{cor}[theorem]{Corollary}

\DeclareMathOperator*{\argmax}{arg\,max}

\DeclareMathOperator*{\maximize}{maximize}
\DeclareMathOperator*{\st}{subject\,\, to \quad}
\newcommand{\norm}[1]{\left|\left|#1\right|\right|}
\newcommand{\neighbors}{\mathcal{N}_i\cup\{i\} }
\newcommand{\myspan}{\text{span}}
\newcommand{\conv}{\mathbf{conv}}

\newcommand{\ybar}{\mathbf{\bar{y}}}
\newcommand{\yreg}{\mathbf{y}}
\newcommand{\vreg}{\mathbf{v}}
\newcommand{\xreg}{\mathbf{x}}
\newcommand{\ureg}{\mathbf{u}}

%% file: contents/abstract.tex
\begin{abstract}

We consider a class of discrete optimization problems that aim to maximize a submodular objective function subject to a distributed partition matroid constraint.  More precisely, we consider a networked scenario in which multiple agents choose actions from local strategy sets with the goal of maximizing a submodular objective function defined over the set of all possible actions.  Given this distributed setting, we develop Constraint-Distributed Continuous Greedy (\texttt{CDCG}), a message passing algorithm that converges to the tight $(1-1/e)$ approximation factor of the optimum global solution using only local computation and communication. It is known that a sequential greedy algorithm can only achieve a $1/2$ multiplicative approximation of the optimal solution for this class of problems in the distributed setting. Our framework relies on lifting the discrete problem to a continuous domain and developing a consensus algorithm that achieves the tight $(1-1/e)$ approximation guarantee of the global discrete solution once a proper rounding scheme is applied.  We also offer empirical results from a multi-agent area coverage problem to show that the proposed method significantly outperforms the state-of-the-art sequential greedy method.

\end{abstract}

%% file: contents/intro.tex
\section{Introduction}

Recently, the need has arisen to design algorithms that distribute decision making among a collection of agents or computing devices.  This need has been motivated by problems from statistics, machine learning and robotics.  More specifically, these problems include:
\begin{itemize}[noitemsep]
    \item \textit{(Density estimation)} What is the best way to estimate a non-parametric density function from a distributed dataset? \citep{hu2007distributed}
    \item \textit{(Non-parametric models)} How should we summarize very large datasets in a distributed manner to facilitate Gaussian process regression? \citep{mirzasoleiman2016distributed}
    \item \textit{(Information acquisition)} How should a team of mobile robots acquire information about an environmental process or reduce uncertainty in a mapping task? \citep{schlotfeldt2018anytime}
\end{itemize}
Research toward solving the problems posed in these applications has resulted in a large body of work on topics such as sensing and coverage \citep{zhong2011distributed, singh2009efficient}, natural language processing \citep{wei2013using}, and learning and statistics \citep{golovin2011adaptive, djolonga2016variational}.   Indeed, inherent to each of these applications is an underlying optimization problem that can be expressed as
\begin{subequations}
\label{eq:intro_statement}
\begin{alignat}{2}
    &\maximize  &&f(\setS) \\
    &\st && \setS\subseteq \Y, \hspace{2pt} \setS\in\I
\end{alignat}
\label{eq:basic-formulation}%
\end{subequations}
where $f$ is a submodular set function (i.e.\ it has a diminishing-returns property), $\Y$ is a finite set of all decision variables, and $\I$ is a family of allowable subsets of $\Y$.  In words, the goal of (\ref{eq:basic-formulation}) is to pick a set $\setS$ from the family of allowable subsets $\mathcal{I}$ that maximizes the submodular set function $f$.  A wide class of relevant objective functions such as  mutual information and weighted coverage are submodular; this has motivated a growing body of work surrounding submodular optimization problems \citep{mokhtari2018decentralized,mirzasoleiman2013distributed,zhou2020distributed,du2020jacobi,adibi2020submodular,chen2020black,xie2019decentralized}.  

Intuitively, it is useful to think of the problem in (\ref{eq:intro_statement}) as a distributed $n$-player game.  In this game, each player or agent has a distinct local strategy set of actions.  The goal of the game is for each agent to choose at most one action from its own strategy set to maximize a problem-specific notion of reward.  Therefore, the problem is \emph{distributed} in the sense that agents can only form a control policy with the actions from their local, distinct strategy sets.  To maximize reward, agents are allowed to communicate with their direct neighbors in a bidirectional communication graph.  In this way, we might think of these agents as robots that collectively aim to solve a coverage problem in an unknown environment by communicating their sensing actions to their nearest neighbors. Throughout this work, we will refer to this multi-agent game example to elucidate our results.

In this paper, our aim is to study problem (1) in a \emph{distributed} setting, which we will formally introduce in Section \ref{sect:prob-statement}; this setting differs considerably from the \emph{centralized} setting, which has been studied thoroughly in past work (see \citet{calinescu2011maximizing}).  Notably, the distributed setting admits a more challenging problem because agents can only communicate locally with respect to a communication graph.  Therefore designing an efficient communication scheme among agents is a concomitant requirement for the distributed setting, whereas in the centralized setting, there is no such desideratum.

\paragraph{Contributions.}  In this paper, we formulate the general case of maximizing a submodular set function subject to a distributed partition matroid constraint in Problem \ref{prob:main_problem}.  We then formulate the continuous relaxation of this problem via the multilinear extension in Problem \ref{prob:cont}.  Both of these problems are formally defined in Section \ref{sect:prob-statement}.  To this end, we study the special case of this optimization problem in which each agent can compute the global objective function and the gradient of the objective function; however we assume that each agent only has access to a local, distinct set of actions.  Considering these constraints, we develop Constraint-Distributed Continuous Greedy (\texttt{CDCG}), a novel algorithm for solving the continuous relaxation of the \emph{distributed} submodular optimization problem that achieves a tight $(1-1/e)$ approximation of the optimal solution, which is known to be the best possible approximation unless $\mathbf{P} = \mathbf{NP}$.  We offer an analysis of the proposed algorithm and prove that it achieves the tight $(1-1/e)$ approximation and that its error term vanishes at a linear rate.  

Previous work on the distributed version of this problem can approximate the optimal solution to within a multiplicative factor of $1/2$ via sequential greedy algorithms \citep{gharesifard2017distributed,corah2018distributed,calinescu2011maximizing}.  Algorithms for different settings, such as the setting of \citep{mokhtari2018decentralized} in which each node has access to a local objective function which is averaged to form a global objective function, can also achieve the $(1-1/e)$ approximation.  Similarly, \citep{calinescu2011maximizing} shows that it is possible to achieve the optimal $(1-1/e)$ approximation in the centralized setting.  However, to the best of our knowledge the \texttt{CDCG} algorithm presented in this paper is the first algorithm that is guaranteed to achieve the $(1-1/e)$ approximation of the optimal solution in this distributed setting.  

%% file: contents/related-work.tex
\section{Related work}

The optimization problem in (\ref{eq:intro_statement}) has previously been studied in settings that differ significantly from the setting studied in this paper.  In particular, \citep{calinescu2011maximizing} addresses this problem in a centralized setting and shows that a centralized algorithm can obtain the tight $(1-1/e)$ approximation of the optimal solution.  In this way, \citep{calinescu2011maximizing} is perhaps the closest to this paper in that both manuscripts introduce algorithms that obtain the tight $(1-1/e)$ guarantee for solving the optimization problem in \eqref{eq:basic-formulation} with respect to a particular setting.  However, the setting of \citep{calinescu2011maximizing} is inherently centralized, whereas our setting is \emph{distributed}. 

Another similar line of work concerns the so-called ``master-worker'' model.  In this framework, agents solve a distributed optimization problem such as \eqref{eq:intro_statement} by exchanging local information with a centralized master node.  However, this setting also differs from the setting studied in this work in that our results assume an entirely distributed setting with no centralized node \citep{mirzasoleiman2013distributed, barbosa2015power}.

Fundamentally, the optimization problem posed in (\ref{eq:intro_statement}) is NP-hard.  However, near-optimal solutions to (\ref{eq:basic-formulation}) can be approximated by greedy algorithms \citep{nemhauser1978analysis, nemhauser1978best}.  In the distributed context, the sequential greedy algorithm (SGA) has been rigorously studied in \citep{gharesifard2017distributed}.  This work poses (\ref{eq:intro_statement}) as a communication problem among agents distributed in an directed acyclic graph (DAG) working to optimize a global objective function.  The authors of \citep{gharesifard2017distributed} offer upper and lower bounds on the performance of SGA based on the clique number of the underlying DAG.  Building on this, \citep{corah2018distributed} analyzes the communication redundancy in such an approach and proposes a distributed planning technique that randomly partitions the agents in the DAG. On the other hand, \citep{grimsman2018impact} extends the work of \citep{gharesifard2017distributed} to a sequential setting in which agents have limited access to the prior decisions of other agents.  Extensions of SGA such as the distributed SGA (DSGA) have also been proposed.  In particular, \citep{corah2017efficient, corah2019distributed} pose (\ref{eq:intro_statement}) as a multi-robot exploration problem and uses DSGA to quantify the suboptimality incurred by redundant sensing information. 

Others have proposed novel algorithms with the goal of avoiding the communication overhead incurred by deploying SGA for a large number of agents.  Instead of explicitly solving (\ref{eq:intro_statement}), many of these algorithms seek to solve a continuous relaxation of this problem \citep{hassani2017gradient,mokhtari2020stochastic}.  This continualization of the problem in (\ref{eq:intro_statement}) was originally introduced in \citep{calinescu2011maximizing}.  In particular, \citep{mokhtari2018decentralized} proposes several gradient ascent-style algorithms for solving a problem akin to (\ref{eq:intro_statement}) in which each agent has access to a local objective function.  Similarly, novel algorithms have been developed for solving problems such as unconstrained submodular maximization \citep{buchbinder2015tight} and submodular maximization with matroid constraints \citep{calinescu2011maximizing,buchbinder2014submodular} by first lifting these problems to the continuous domain.

Another notable direction in solving problem \eqref{eq:basic-formulation} has been to define an auxiliary or surrogate function in place of the original submodular objective.  For instance, \citep{clark2015scalable} introduces a distributed algorithm for maximizing a submodular auxiliary function subject to matroid constraints that obtains the $(1-1/e)$ optimal approximation.  This approach of defining surrogate functions in place of the submodular objective differs significantly from our approach.

%% file: contents/prelims.tex
\section{Preliminaries}

In this section, we review the notation used throughout this paper and state definitions that are necessary for the problem formulations in Section \ref{sect:prob-statement}. 

\paragraph{Notation.}
Throughout this paper, lowercase bold-face (e.g.\ $\vreg$) will denote a vector, while uppercase bold-face (e.g.\ $\W$) will denote a matrix.  The $i^{\text{th}}$ component of a vector $\vreg$ will be denoted $v_i$; the element in the $i^{\text{th}}$ row of the $j^{\text{th}}$ column of a matrix $\W$ will be denoted by $w_{ij}$.  The inner product between two vectors $\xreg$ and $\yreg$ will be denoted by $\langle \xreg, \yreg \rangle$ and the Euclidean norm of a vector $\vreg$ will be denoted by $\norm{\vreg}$.  Given two vectors $\xreg$ and $\yreg$, we define $\xreg\vee\yreg = \max(\xreg, \yreg)$ as the (vector-valued) component-wise maximum between $\xreg$ and $\yreg$; similarly, $\xreg\wedge\yreg = \min(\xreg, \yreg)$ will denote the component-wise minimum between $\xreg$ and $\yreg$.  We will use the notation $\0_n$ to denote an $n$-dimensional vector in which each component is zero; similarly $\1_n$ will denote an $n$-dimensional vector in which each component is one.  Calligraphic fonts will denote sets (e.g.\ $\Y$).  Given a set $\Y$, $|\Y|$ will denote the cardinality of $\Y$, while $2^{\Y}$ will denote the power set of $\Y$.  $\indicator_{\Y}:\Y\mapsto\{0,1\}$ will represent the indicator function for the set $\Y$.  That is, $\indicator_{\Y}$ is the function that takes value one if its argument is an element of $\Y$ and takes value zero otherwise.  Finally, $\varnothing$ will denote the null set.


\paragraph{Background and relevant definitions.}
Let $\Y$ be a finite set and let $f:2^{\Y}\mapsto\R_+$ be a set function mapping subsets of $\Y$ to the nonnegative real line.  In this setting, $\Y$ is commonly referred to as the \textit{ground set}.  The function $f$ is called \textit{submodular} if for every $\A, \B\subseteq \Y$,
\begin{align*}
    f(\A\cap \B) + f(\A\cup \B) \leq f(\A) + f(\B).
\end{align*}
In essence, submodularity amounts to $f$ having a so-called diminishing-returns property, meaning that the incremental value of adding a single element to the argument of $f$ 
is no less than that of
adding the same element to a superset of the argument.
To illustrate this, we will slightly overburden our notation by defining
\begin{align*}
    f(\xreg|\A) := f(\A \cup \{\xreg\}) - f(\A)
\end{align*}
as the \textit{marginal reward} of $x$ given $\A$.  This gives rise to an equivalent definition of submodularity.  In particular, $f$ is said to be submodular if for every $\A\subseteq \B\subseteq\Y$ and $\forall \xreg\in\Y\backslash\B$,
\begin{align*}
    f(\xreg|\B) \leq f(\xreg|\A).
\end{align*}
Throughout this paper, we will consider submodular functions that are also \textit{monotone}, meaning that for every $\A\subseteq \B\subseteq \Y$, $f(\A)\leq f(\B)$, and \textit{normalized}, meaning that $f(\varnothing) = 0$.

In practice, one often encounters a constraint on the allowable subsets of the ground set $\Y$ when maximizing a submodular objective function.  Concretely, if $\I$ is a nonempty family of allowable subsets of the ground set $\Y$, then the tuple $(\Y, \I)$ is a \textit{matroid} if the following criteria are satisfied:
\begin{enumerate}
    \item[(1)] (\textit{Heredity}) For any $\A\subset \B\subset \Y$, if $\B\in\I$, then $\A\in\I$.
    \item[(2)] (\textit{Augmentation}) For any $\A, \B\in\I$, if $|\A|<|\B|$, then $\exists$ $\xreg\in\B\backslash \A$ such that $\A \cup\{\xreg\} \in\I$.
\end{enumerate}
Furthermore, if $\Y$ is partitioned into $n$ disjoint sets $\Y_1, \dots, \Y_n$, then the tuple $(\Y, \I)$ is a \textit{partition matroid} if there exists positive integers $\alpha_1, \dots, \alpha_n$ such that
\begin{align*}
    \I \equiv \{\A \hspace{3pt} : \hspace{3pt} \A \subseteq \Y, \hspace{1pt} |\A \cap \Y_i| \leq \alpha_i \hspace{3pt} \text{for each } i = 1, \dots, n\}.
\end{align*}
Partition matroids are particularly useful when defining the constraints of a distributed optimization problem because they can be used to describe a setting in which a ground set $\Y$ of all possible actions is written as the product of disjoint local action spaces $\Y_i$.

The notion of submodularity can be extended to the continuous domain \citep{wolsey1982analysis}.  Consider a set $\X = \prod_{i=1}^n \X_i$, where $\X_i$ is a compact subset of $\R_+$ for each index $i \in \{1, \dots, n\}$.  We call a continuous function $F:\X\rightarrow\R_+$ \textit{submodular} if for all $\xreg, \yreg \in \X$,
\begin{align*}
    F(\xreg \vee \yreg) + F(\xreg \wedge \yreg) \leq F(\xreg) + F(\yreg).
\end{align*}
As in the discrete case, we say that a continuous function $F$ is monotone if $\forall\xreg, \yreg\in\X$, $\xreg\preceq \yreg$ implies that $F(\xreg) \leq F(\yreg)$.  Furthermore, if $F$ is differentiable, we say that $F$ is $DR$-submodular, where $DR$ stands for ``diminishing-returns,'' if the gradients are \textit{antitone}.  That is, $\forall\xreg,\yreg \in\X$, $F$ is $DR$-submodular if $\xreg\preceq\yreg$ implies that $\nabla F(\xreg) \succeq F(\yreg)$. 

%% file: contents/prob-statement.tex

\section{Problem Statement}
\label{sect:prob-statement}

In this section, we formulate the main problem of this paper: maximizing submodular set functions subject to distributed partition matroid constraints.

\begin{problem}[\textbf{Submodular Maximization Subject to a Distributed Partition Matroid Constraint}]\label{prob:main_problem}
Consider a collection of $n$ agents that form the set $\agents=\{1, \dots, n\}$.  Let $f:2^{\Y} \mapsto \R_{+}$ be a normalized and monotone submodular set function and let $\Y_1, \dots, \Y_n$ be a pairwise disjoint partition of a finite ground set $\Y$, wherein each agent $i\in\agents$ can only choose actions from its local strategy set $\Y_i$.  Furthermore, consider the  partition matroid $(\Y, \I)$, where
\begin{align}
    \I := \{\setS\subseteq \Y \hspace{3pt} : \hspace{3pt}  |\Y_i\cap\setS|\leq 1 \hspace{3pt} {\normalfont \text{for } } i = 1, \dots, n \}
    \label{eq:main_problem_matroid_constraint}.
\end{align}
The problem of submodular maximization subject to a distributed partition matroid constraint is to maximize $f$ by selecting a set $\setS\subseteq\Y$ from the family of allowable subsets so that $\setS \in\I$.  Formally:
\vspace{-1.5em}
\begin{subequations}
\label{eq:main_problem}
\begin{alignat}{2}
    &\maximize &&f(\setS) 
    \label{eq:main_problem_objective} \\
    &\st && \setS \in \I
    \label{eq:main_problem_constrain}
\end{alignat}%
\end{subequations}
\end{problem} 
In effect, the distributed partition matroid constraint in Problem \ref{prob:main_problem} enforces that each agent $i\in\agents$ can choose at most one action from its local strategy set $\Y_i$.  Note that in this setting, each agent can only choose actions from its own local strategy set.  Therefore, this problem is distributed in the sense that agents can only determine the actions taken by other agents by directly communicating with one another.

\subsection{Sequential greedy algorithm}

It is well known that the sequential greedy algorithm (SGA), in which each agent $i\in\agents$ chooses an action sequentially based on 
\begin{align}
    \yreg_i = \argmax_{\yreg \in \Y_i} f(\yreg | \setS_{i-1})
\end{align}
where $\setS_{i-1} = \{\yreg_1, \dots, \yreg_{i-1}\}$, approximates the optimal solution to within a multiplicative factor of $1/2$ \citep{gharesifard2017distributed}.  The drawbacks of this algorithm are twofold.  Firstly, as we will show, our algorithm achieves the tight $(1-1/e)$ approximation of the optimal solution, which is known to be the best possible approximation unless $\mathbf{P} = \mathbf{NP}$.  Secondly, as its name suggests, SGA is sequential in nature and therefore it scales very poorly in the number of agents.  That is, each agent must wait for each of the previous agents to compute their contribution to the optimal set $\setS^*$.  Notably, our algorithm does not suffer from this sequential dependence.

\subsection{Continuous Extension of Problem \ref{prob:main_problem}}
\label{subsect:cont_ext}

Sequential algorithms such as SGA can only achieve a $1/2$ approximation of the optimal solution.  To achieve the best possible $(1-1/e)$ approximation of the optimal solution, it is necessary to extend Problem \ref{prob:main_problem} to the continuous domain via the so-called \textit{multilinear extension} of the submodular objective function $f$ \citep{nemhauser1978analysis}.  Thus, the method we use in this work to achieve the tight $(1-1/e)$ approximation relies on the continualization of Problem \ref{prob:main_problem}.  Importantly, it has been shown that Problem \ref{prob:main_problem} and the optimization problem engendered by lifting Problem 1 to the continuous domain via this multilinear extension yield the same solution \citep{calinescu2011maximizing}.  Furthermore, by applying proper rounding techniques, such as those described in Section 5.1 of \citep{mokhtari2018decentralized} and in \citep{calinescu2011maximizing} and \citep{chekuri2014submodular} to the continuous relaxation of Problem \ref{prob:main_problem}, one can obtain the tight $(1-1/e)$ approximation for Problem \ref{prob:main_problem}.  Therefore, our approach in this paper will be to lift Problem \ref{prob:main_problem} to the continuous domain.  We formulate this problem in the following way: 

\begin{problem}[\textbf{Continuous Extension of Problem \ref{prob:main_problem}}] \label{prob:cont}
Consider the conditions of Problem \ref{prob:main_problem}.  Define the $DR$-submodular continuous multilinear extension $F:\X\mapsto\R_+$ of the objective function $f$ in Problem \ref{prob:main_problem} by 
\begin{align}
    F(\yreg) := \sum_{\setS \subseteq \Y} f(\setS) \prod_{i\in\setS} y_i \prod_{j\not\in\setS} (1-y_j)
    \label{eq:cont_ext_obj_func}
\end{align}
and let $\setP\subseteq\X$ be the matroid polytope $\setP := \conv\{\indicator_{\setS} \hspace{3pt} : \hspace{3pt} \setS\in\I \}$ where $\I$ is the family of sets defined in (\ref{eq:main_problem_matroid_constraint}).  The continuous relaxation of Problem \ref{prob:main_problem} is formally defined by
\begin{subequations}
\label{eq:cont_ext_problem}
\begin{alignat}{2}
    &\maximize && F(\yreg) \\
    &\st && \yreg\in \setP
    \label{eq:cont_prob_constraint}
\end{alignat}
\end{subequations}
\end{problem}
Observe that Problem \ref{prob:cont} is \emph{distributed} in the sense that each agent $i\in\agents$ is associated with its own distinct continuous strategy space $\setP_i$.  Formally, the set $\setP_i$ is defined as
\begin{align}
    \setP_i := \conv\{\indicator_{\setS} \hspace{3pt} : \hspace{3pt} S\subseteq \I_i\}
    \label{eq:def_p_i}
\end{align}
where $\I_i := \{ \setS\subseteq \Y \hspace{3pt} : \hspace{3pt} |\Y_i \cap \setS| \leq 1\}$.  In this way, $\setP = \cap_{i=1}^n\setP_i$.  In this way, the sets $\setP_i$ play similar roles in Problem \ref{prob:cont} as the sets $\Y_i$ do in Problem \ref{prob:main_problem}.

Note that Problem \ref{prob:cont} is nonconvex, and therefore cannot be solved by classical convex solvers or algorithms.  Further, we assume that each agent $i\in\agents$ can compute the multilinear extension $F$ of the submodular objective function $f$ in (\ref{eq:main_problem_objective}) and the gradient of $F$.  

%% file: contents/algorithm.tex

\section{Constraint-Distributed Continuous Greedy}
\label{sect:algo-desc}

In this section, we present Constraint-Distributed Continuous Greedy (\texttt{CDCG}), a decentralized algorithm for solving Problem \ref{prob:cont}.  The pseudo-code of \texttt{CDCG} is described in Algorithm \ref{algo:cdcg}.  At a high level, this algorithm involves updating each agent's local decision variable based on the aggregated belief of a small group of other agents about the best control policy.  In essence, inter-agent communication within small groups of agents facilitates local decision making.

For clarity, we introduce a simple framework for the inter-agent communication structure.  In \texttt{CDCG}, agents $i\in\agents = \{1, \dots, n\}$ share their decision variables $\yreg_i$ with a small subset of \textit{local} agents in $\agents$.  To encode the notion of locality, suppose that each agent $i\in\agents$ is a node in a bidirectional \textit{communication graph} $\G = (\agents, \E)$ in which $\E$ denotes the set of edges.  Given this structure, we assume that each agent $i\in\agents$ can only communicate its decision variable $\yreg_i$ with its direct neighbors in $\G$.  Let us denote the neighbor set of agent $i\in\agents$ by $\agents_i$.  Then the set of edges $\E$ can be written $\{(i,j) \hspace{3pt} : \hspace{3pt} j\in\agents_i\}$.  We adopt this notation for the remainder of this paper.

\subsection{Intuition for the \texttt{CDCG} algorithm}

The goal of \texttt{CDCG} at a given node $i\in\agents$ is to learn the local decision variable $\yreg_i$.  \texttt{CDCG} is run at each node in $i\in\agents$ to assemble the collection $\{\yreg_1^T, \dots, \yreg_n^T\}$ where $T$ is a given positive integer; this collection represents an approximate solution to Problem \ref{prob:cont} and guarantees that each agent contributes at most one element to the solution.  Then, by applying proper rounding techniques to each element of the collection such as those discussed in \citep{mokhtari2018decentralized,calinescu2011maximizing,chekuri2014submodular}, we obtain a solution to Problem \ref{prob:main_problem}.  In the proceeding sections, we show that this solution achieves the tight $(1-1/e)$ approximation of the optimal solution. 

In the analysis of \texttt{CDCG}, we add the superscript $t$ to the vectors $\vreg_i^t$ and $\yreg_i^t$ defined in Algorithm \ref{algo:cdcg}.  This superscript denotes the iteration number so that $\yreg_i^t$ and $\vreg_i^t$ represent the values of the local variables $\yreg_i$ and $\vreg_i$ at iteration $t\in\{1, \dots, T\}$ respectively.

\subsection{Description of the steps for \texttt{CDCG} (Algorithm \ref{algo:cdcg})}

From the perspective of node $i\in\agents$, \texttt{CDCG} takes two arguments: nonnegative weights $w_{ij}$ for each $j\in\neighbors$ and a positive integer $T$.  The weights $w_{ij}$ correspond to the $i^{\text{th}}$ row in a doubly-stochastic weight matrix $\W$ and $T$ is the number of iterations for which the algorithm will run.  The weight matrix $\W$ is a design parameter of the problem and must fulfill a number of technical requirements that are fully described in Appendix A. 
Before any computation, the local decision variable $\yreg_i$ is initialized to the zero vector.

Computation proceeds in $T$ rounds.  In each round, the first step is to calculate the gradient of the multilinear extension function $F$ evaluated at the local decision variable $\yreg_i^{t-1}$ from the previous iteration.  Thus, in line 3 of Algorithm \ref{algo:cdcg}, we calculate the ascent direction $\vreg_i^t$ at iteration $t$ in the following way:
\begin{align*}
    \vreg_i^t = \argmax_{\xreg \in\setP_i \cap \setC_i} \left\langle \nabla F(\yreg_i^{t-1}), \xreg \right\rangle.
\end{align*}
Intuitively, one can think of $\vreg_i^t$ as the vector from the set $\setP_i\cap\setC_i$ that is most aligned with $\nabla F(\yreg_i^{t-1})$.  To define the set $\setC_i$, first define the set $\setJ_i$ as the set of indices of the elements in $\Y$ that correspond to elements in $\Y_i$.  Then
\begin{align}
    \setC_i := \left\{ \xreg \in\R_+^{|\Y|} \hspace{3pt} : \hspace{3pt} x_j = 0 \quad\forall j\not\in\setJ_i \right\}.
    \label{eq:def_c_i}
\end{align}
Using this notation, we can equivalently define $\setP_i = \{ \xreg\in\R_+^{|\Y|} \hspace{3pt} : \hspace{3pt} \sum_{j\in\setJ_i} x_j \leq 1 \}$.  Next, in line 4 of Algorithm \ref{algo:cdcg}, $\yreg_i$ is updated by setting
\begin{align*}
    \yreg_i^t = \sum_{j\in\neighbors} w_{ij}\yreg_j^{t-1} + \frac{n}{T} \vreg_i^t.
\end{align*}
In this way, the governing principle is to collaboratively accumulate the local belief about the optimal decision $\yreg_i^{t-1}$ and to then move in the approximate direction of steepest ascent from this point.

After $T$ rounds of computation at each node $i\in\agents$, we obtain a local decision variable $\yreg_i^T$ at each node.  By applying proper rounding techniques, we obtain a decision variable for each agent $i\in\agents$.  Rounding in a decentralized manner is discussed in Section 5.1 of \citep{mokhtari2018decentralized}.  The rounding techniques of \citep{mokhtari2018decentralized} build on ``pipage rounding'' \citep{calinescu2011maximizing} and ``swap rounding'' \citep{chekuri2014submodular}, which are both centralized rounding techniques.  The collection of these decision variables form the set $\setS^*$, which represents our solution to Problem~\ref{prob:main_problem}.

\begin{algorithm}[t]
\caption{Constraint-Distributed Continuous Greedy (\texttt{CDCG}) at node $i$}
\textbf{Require:} Weights $w_{ij}$ for each neighbor $j\in\neighbors$ and number of rounds $T\in\Z_{++}$\\
\textbf{Returns:} Local solution $\mathbf{y}_i^\star$ for node $i\in\mathcal{N}$ to Problem \ref{prob:main_problem}
\begin{algorithmic}[1]
    \STATE Initialize local vectors $\mathbf{y}_i^0 = \0_{|\Y|}$
    \FOR{$t = 1, 2, \dots, T$}
        \STATE \textbullet \quad Calculate an ascent direction for the multilinear extension function $F$ via:
        \vspace{-0.5em}
        \begin{align*}
            \mathbf{v}_i^t \gets \mathop{\argmax}\limits_{\xreg\in \setP_i \cap \setC_i} \left\langle \nabla F\left(\yreg_i^{t-1}\right), \xreg \right\rangle
        \end{align*}
        \vspace{-0.5em}
        \STATE \textbullet \quad Update the local variable $\mathbf{y}_i^t$ using the ascent direction $\mathbf{v}_i^t$ via:
        \vspace{-0.5em}
        \begin{align*}
            \mathbf{y}_i^t \gets \displaystyle\sum_{j\in\neighbors} w_{ij}\mathbf{y}_j^{t-1} + \frac{n}{T} \mathbf{v}_i^t
        \end{align*}
        \vspace{-1.0em}
    \ENDFOR
    \STATE $\mathbf{y}_i^\star \gets$ \texttt{Round}$\left(\mathbf{y}_i^T\right)$
\end{algorithmic}\label{algo:cdcg}
\end{algorithm}

%% file: contents/convergence.tex
\section{Convergence Analysis}

The main result in this paper is to show that in the distributed setting of Problem \ref{prob:cont}, \texttt{CDCG} achieves a tight $(1-1/e)$ multiplicative approximation of the optimal solution.  The following theorem summarizes this result.  

\begin{theorem}
\label{thm:main_result}
Consider the \texttt{CDCG} algorithm described in Algorithm \ref{algo:cdcg}.  Let $\yreg^*$ denote the global maximizer of the optimization problem defined in Problem \ref{prob:cont}, and assume that a positive integer $T$ and a doubly-stochastic weight matrix $\W$ are given.  Then provided that the assumptions outlined in Appendix A hold, for all nodes $i\in\mathcal{N}$, the local variables $\yreg_i^T$ obtained after $T$ iterations satisfy
\begin{align}
    F(\yreg_i^T) &\geq \left( 1 - \frac{1}{e} \right) F(\yreg^*) - \left[\frac{LD^2}{2T} + \frac{LD^2(n^2 + n^{5/2}) + n^{5/2}DG }{T(1-\beta) }\right] 
    \label{thm:main_result_statement}
\end{align}
where $D$, $G$, $L$, and $\beta$ are problem-dependent constants that are formally defined in Appendices A and B.  
\end{theorem}

\noindent Succinctly, Theorem \ref{thm:main_result} means that the sequence of local iterates generated by \texttt{CDCG} achieves the optimal approximation ratio $(1-1/e)$ and that the error term vanishes at a linear rate of $\bigO(1/T)$.  That is,
\begin{align*}
    F(\yreg_i^T) \geq \left(1-\frac{1}{e}\right) F(\yreg^*) - \bigO\left(\frac{1}{T} \right),
\end{align*}
which implies that each agent reaches an objective value larger than $(1-1/e-\epsilon)\yreg^*$ after $\bigO(1/\epsilon)$ rounds of communication.  Previous work can only guarantee an objective value of $(1/2)\yreg^*$ \citep{gharesifard2017distributed}.  We provide the proof of this theorem and supporting lemmas in Appendices B and C.  

%% file: contents/simulation.tex

\section{Simulation Results}

\begin{figure}[t]
    \centering
    \includegraphics[width=0.95\textwidth]{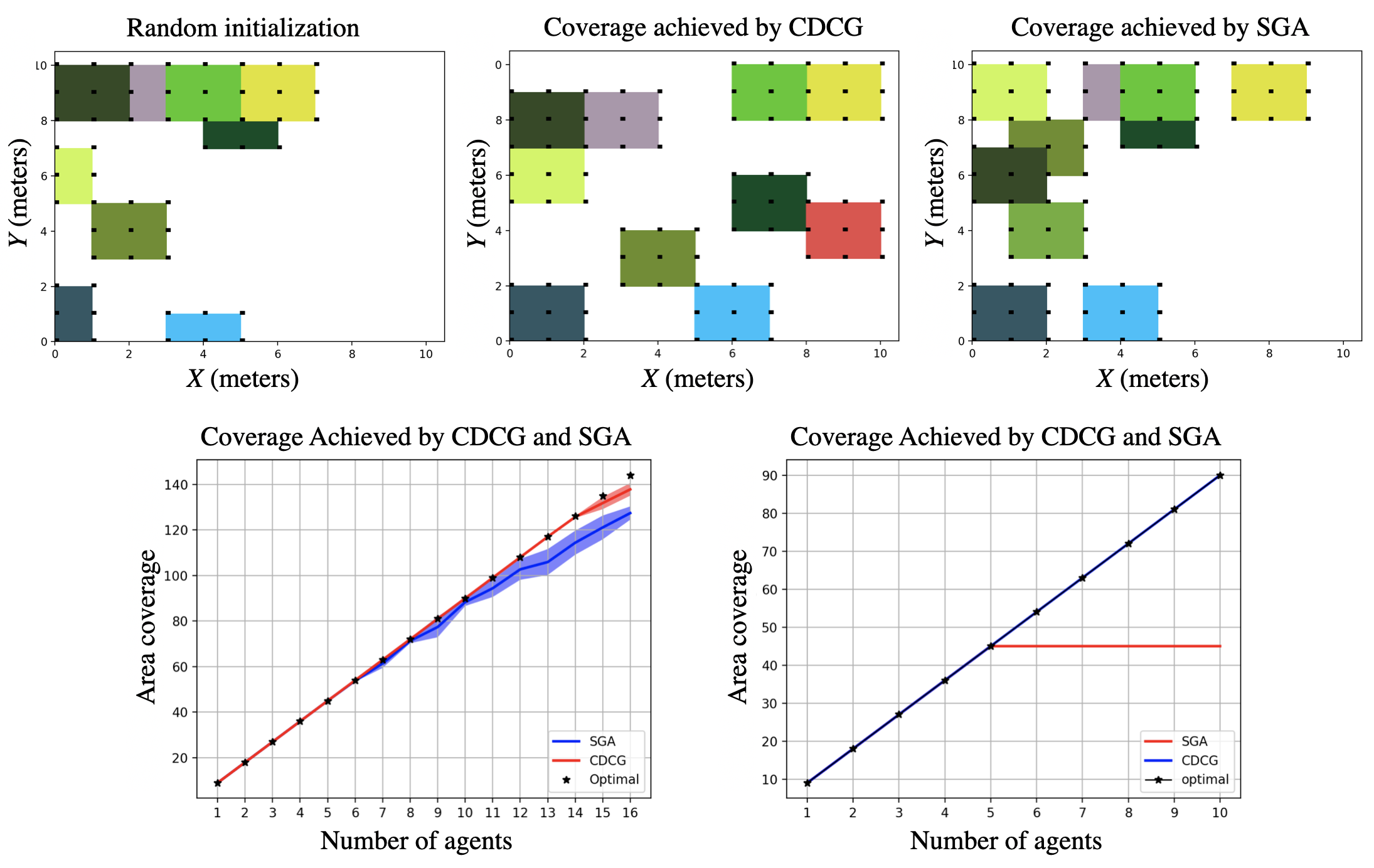}
    \vspace{-1.0em}
    \caption{\textbf{Area coverage simulation results for \texttt{CDCG} and \texttt{SGA}.}  (Top left) Random initialization of $n=10$ agents in a $10\times 10$ grid.  (Top middle \& right) Coverage achieved by \texttt{CDCG} (top middle) and \texttt{SGA} (top right) from the random initialization shown in the top left panel.  (Bottom left) Comparison of the mean coverage achieved by \texttt{CDCG} and \texttt{SGA} averaged over 10 random initializations.  (Bottom right) Comparison of the coverage achieved by \texttt{CDCG} and \texttt{SGA} for a setting in which each agent's starting point is the center of the grid.}
    \label{fig:coverage}
\end{figure}

To evaluate the proposed algorithm, we consider a multi-agent area coverage problem.  In this setting, each agent $i\in\agents$ is constrained to move in a two-dimensional grid.  We assume that each agent has a finite radius $r$ so that it can observe those grid points that lie with a square with sidelength $2r+1$.  The objective is for the agents to collectively maximize the cardinality of the union of their observation sets of grid points.  In other words, given an initial configuration, the problem is to choose an action for each agent that maximizes the overall coverage of the grid.  The top three panels of Figure \ref{fig:coverage} show various configurations of agents in this two-dimensional grid.

Consider an initial configuration of $n$ agents in states $\yreg_{i} \in \mathbb{Z}^2$ for $i\in\{1, \dots, n\}$ with the dynamic constraint $\yreg_{i}^{t+1} = \yreg_i^t + \ureg_i^t$, where $\ureg_i^t$ is a control input from a discrete set 
\begin{align*}
    \mathcal{U} = \{(0, 1), (0, -1), (-1, 0), (1, 0), (0, 0)\}.
\end{align*}
Elements from this set represent the admissible actions for each agent in the two-dimensional grid.

In our simulation, we compared the performance of \texttt{SGA} against \texttt{CDCG} on the coverage task posed above for a variable number of agents.  For simplicity, we assumed that the underlying communication graph $\mathcal{G}$ used in \texttt{CDCG} was fully connected and that each value in the weight matrix $\W$ was $1/n$.  A random initialization for each agent's position and the coverage achieved by \texttt{CDCG} and \texttt{SGA} are shown in the top three panels of Figure \ref{fig:coverage} respectively.  We compared the performance of these algorithms across ten random initializations of starting locations for the agents; the mean performance of each algorithm and the respective standard deviations are shown in the bottom left panel of Figure \ref{fig:coverage}.  In each trial, we ran both algorithms 50 times, each of which produced a control input $\ureg_i$ for each agent.  For each initialization, we ran \texttt{CDCG} for $T=100$ iterations.  Note that as the number of agents increases, \texttt{CDCG} is optimal or near optimal in each case; however for larger than eight agents, the performance of \texttt{SGA} begins to fall away from the optimal.

We also compared the coverages achieved by \texttt{CDCG} and \texttt{SGA} for a setting in which each agent's starting position is the center of the grid.  The results of this experiment are shown in the bottom right panel of Figure \ref{fig:coverage}.  In this plot, we averaged the performance over 15 independent trials; in each trial, we ran $\texttt{CDCG}$ for $T=100$ iterations.  Interestingly, \texttt{SGA} converges to a local maximum in this problem, whereas \texttt{CDCG} achieves the optimal value.

%% file: contents/conclusion.tex
\section{Conclusion}

In this work, we described an approach for achieving the optimal approximation to a class of submodular optimization problems subject to a \emph{distributed} partition matroid constraint.  The algorithm we proposed outperforms the sequential greedy algorithm in two senses: (1) \texttt{CDCG} achieves the tight $(1-1/e)$ approximation for the optimal solution whereas SGA can only achieve a $1/2$ approximation; and (2) \texttt{CDCG} imposes a limited communication structure on this problem, which allows for significant gains via parallelization.  We showed empirically via an area coverage simulation with multiple agents that \texttt{CDCG} outperforms the greedy algorithm.


%% file: appendix/assumptions.tex
\section*{Appendix A: Assumptions for Theorem \ref{thm:main_result}}

Consider the continuous relaxation of Problem \ref{prob:main_problem} that was described in Section \ref{subsect:cont_ext}.  We assume that the Euclidean distance between elements of the convex set $\setP$ are uniformly bounded, i.e. that
\begin{align}
    \norm{\xreg - \yreg} \leq D
    \label{assume:unif-bounded} \qquad \forall  \xreg,\yreg \in\setP.
\end{align}
This is a trivial consequence of the multilinear extension $F$, since $\mathcal{P}$ is contained in the unit cube.  Furthermore, we assume that the gradient of the multilinear extension $F$ of the objective function $f$ in Problem \ref{prob:main_problem} is $L$-Lipschitz continuous, i.e. that
\begin{align}
    \norm{\nabla F(\xreg) - \nabla F(\yreg)} \leq L\norm{\xreg - \yreg} \qquad \forall \xreg,\yreg \in\setP 
    \label{assume:l-lipschitz-grad}
\end{align}
so that $\norm{\nabla F(\xreg) - \nabla F(\yreg)} \leq LD$ $\forall \xreg,\yreg \in\setP$ by (\ref{assume:unif-bounded}).  Again, this is not a limiting assumption, because the domain of $F$ is compact, which implies the Lipschitzness of $F$. Also, we assume that the norm of the gradient of $F$ is bounded over $\setP$, i.e. that
\begin{align}
    \norm{\nabla F(\xreg)} \leq G \qquad\forall \xreg\in\setP
    \label{assume:bounded-grad},
\end{align}
which again follows from the compactness of the domain of $F$.  It is then easy to show that (\ref{assume:bounded-grad}) and the multivariable mean value theorem imply that $F$ is $G$-Lipschitz continuous over $\setP$.  Note that in this case, since $F$ is the multilinear extension of $f$, assumptions (\ref{assume:unif-bounded}), (\ref{assume:l-lipschitz-grad}), and (\ref{assume:bounded-grad}) all hold.  Moreover, the constants $L$, $D$, and $G$ all depend on the maximum singleton value of $f$.  For further justification, see \citep{hassani2017gradient,mokhtari2018decentralized}.
Finally, it will be prudent to mention that for the multilinear extension $F$ of any monotone and submodular function $f$, it holds that $F(\0) \geq 0$ and
\begin{align}
    \left\langle\nabla F(\ybar^t), \yreg^* \right\rangle \geq F(\yreg^*) - F(\ybar^t)
    \label{assume:concave-obj}
\end{align}
For justification, see \citep{calinescu2011maximizing}.

Now consider the communication framework described in Section \ref{sect:algo-desc} and the weight matrix $\W$.  This matrix is a parameter that is designed to match the criteria and setting of a given application.  We assume that the weights used in \texttt{CDCG} are nonnegative so that $w_{ij} \geq 0$ $\forall i,j\in\agents$; furthermore, if node $j\not\in\agents_i$, then $w_{ij} = 0$.  Also, we assume that the weight matrix $\W$ is doubly stochastic and symmetric, and that $\null(\identity - \W) = \myspan(\1_n)$.  The assumptions made about $\W$ are similar to those described in \citep{mokhtari2018decentralized}.  

Lastly, consider that past work has studied the case in which the objective function is distributed \citep{mokhtari2018decentralized}.  However, our setting is one in which the problem is distributed in the constraints rather than the objective.  Therefore, we assume that each agent has access to an oracle for computing the objective submodular function $f$.

%% file: appendix/prelim-lemmas.tex
\section*{Appendix B: Preliminary Lemmas}
\label{app:proofs}

In this appendix, we offer proofs of lemmas that support the proof of Theorem \ref{thm:main_result}.  We note that the proofs for Lemmas \ref{lemma:convergence} and \ref{lemma:consensus} are similar to those that originally appeared in \citep{mokhtari2018decentralized}, and where relevant, pieces of these arguments have been reproduced for completeness.

In general, the goal of Lemma \ref{lemma:convergence} is to show that the local decision variable $\yreg_i$ for each agent $i\in\agents$ converges to the mean $\ybar = \frac{1}{n}\sum_{i\in\agents} \yreg_i$.  Then, in Lemma \ref{lemma:consensus}, we show that these means are Cauchy, meaning that for a sufficiently large number of iterations $T$, the distance between $\ybar^t$ and $\ybar^{t+1}$ becomes arbitrarily small.  Together, Lemma \ref{lemma:convergence} and Lemma \ref{lemma:consensus} establish that for a sufficiently large number of iterations, the set of nodes come to a consensus for the optimal decision.  Lemmas \ref{lemma:tech_lemma_argmax} and Lemma \ref{lemma:tech_lemma_bound} are technical results used in the proof of Theorem \ref{thm:main_result}.


\begin{lemma}
For any iteration $t\leq T$ where $T\in\Z_{++}$, it follows that the Euclidean distance between the local variable $\yreg_i^t$ at node $i\in\mathcal{N}$ and the mean of the local variables $\ybar^t$ can be bounded by
\begin{align*}
    \norm{\yreg_i^t - \ybar^t} \leq \frac{n^{3/2}D}{T(1-\beta)}
\end{align*}
where $\beta$ is the magnitude of the eigenvalue of $\W$ that among all eigenvalues in $\sigma(\W)$ has the second largest magnitude.

\end{lemma}\label{lemma:convergence}

\begin{proof}
Define $\yregcon := \begin{bmatrix} \yreg_1; \dots; \yreg_n \end{bmatrix} \in\R^{np}$ and $\vregcon := \begin{bmatrix} \vreg_1; \dots; \vreg_n \end{bmatrix} \in\R^{np}$ as the concatenations of the local variables $\yreg_i^t$ and descent directions $\vreg_i$ in \texttt{CDCG}.  The update rule in step 2 in Algorithm \ref{algo:cdcg} leads to the expression
\begin{align}
    \yregcon^t = \frac{n}{T} \sum_{s=0}^{t-1} \left( \W \otimes \identity \right)^{t-1-s} \vregcon^s
    \label{eq:update_rule_kron}
\end{align}
Next, if we premultiply both sides of (\ref{eq:update_rule_kron}) by the matrix $(\frac{\1_n \1_n^{\dagger}}{n} \otimes \identity)$, which is the Kronecker product of the matrices $\frac{\1_n \1_n^{\dagger}}{n} \in\R^{n\times n}$ and $\identity\in\R^{p\times p}$, we obtain
\begin{align}
    \left( \frac{\1_n\1_n^{\dagger}}{n} \otimes \identity \right) \yregcon^t = \frac{n}{T} \sum_{s=0}^{t-1} \left[ \left( \frac{\1_n\1_n^{\dagger}}{n} \W^{t-1-s} \right) \otimes \identity \right] \vregcon^s.
    \label{eq:conv-lem-eq-2}
\end{align}
The left hand side of (\ref{eq:conv-lem-eq-2}) can be simplified to 
\begin{align}
    \left(\frac{\1_n\1_n^{\dagger}}{n} \otimes \identity \right) \yregcon^t = \ybarcon^t
    \label{eq:conv-lem-eq-3}
\end{align}
where $\yregcon^t = \begin{bmatrix} \ybar^t; \dots; \ybar^t \end{bmatrix}$.  Combining (\ref{eq:conv-lem-eq-3}) and the equality $\1_n\1_n^{\dagger}\W = \1_n\1_n^{\dagger}$, we can write (\ref{eq:conv-lem-eq-2}) as
\begin{align}
    \ybarcon^t = \frac{n}{T} \sum_{s=0}^{t-1} \left( \frac{\1_n\1_n^{\dagger}}{n}  \otimes \identity \right) \vregcon^s.
    \label{eq:conv-lem-eq-4}
\end{align}
Using the expressions in (\ref{eq:update_rule_kron}) and (\ref{eq:conv-lem-eq-4}), we can derive an upper bound on the difference $\norm{\yregcon^t - \ybarcon^t}$ by
\begin{align}
    \norm{\yregcon^t - \ybarcon^t} 
    &= \frac{n}{T} \norm{ \sum_{s=0}^{t-1} \left[ \left( \W^{t-1-s} - \frac{\1_n\1_n^{\dagger}}{n}  \right) \otimes \identity \right] \vregcon^s } \notag \\
    &\leq \frac{n}{T} \sum_{s=0}^{t-1} \norm{ \mathbf{W}^{t-1-s} - \frac{\1_n\1_n^{\dagger}}{n} } \cdot \norm{\vregcon^s} \notag \\
    &\leq \frac{nD}{T}  \norm{ \W^{t-1-s} - \frac{\1_n\1_n^{\dagger}}{n} },
    \label{eq:conv-lem-eq-5}
\end{align}
where the first inequality follows from the Cauchy-Schwartz inequality and the fact that the norm of a matrix does not change if we Kronecker it by the identity matrix.  The second inequality holds because $\norm{\vregcon^t} \leq D$.  Note that the eigenvectors of the matrices $\W$ and $\W^{t-1-s}$ are the same for all $s=0, \dots, t-1$.  Therefore, the largest eigenvalue of $\W^{t-1-s}$ is 1 with eigenvector $\1_n$ and the second largest magnitude of the eigenvalues is $\beta^{t-1-s}$, where $\beta$ is the second largest magnitude of the eigenvalues of $\mathbf{W}$.  Also note that because $\1_n$ is an eigenvector of $\W^{t-1-s}$, it follows that all of the other eigenvectors of $\W^{t-1-s}$ are orthogonal to $\1_n$ since $\W$ is symmetric.  Hence we can bound the norm $\norm{\W^{t-1-s} - (\1_n\1_n^{\dagger})/n}$ by $\beta^{t-1-s}$.  Applying this substitution to the right hand side of (\ref{eq:conv-lem-eq-5}) yields
\begin{align}
    \norm{ \yregcon^t - \ybarcon^t } &\leq \frac{nD}{T} \sum_{s=0}^{t-1} \beta^{t-1-s} \leq \frac{nD}{T(1-\beta)}
    \label{eq:conv_lemma_final_bound}.
\end{align}
Since $\norm{\yregcon^t - \ybarcon^t}^2 = \sum_{i=1}^n \norm{\yreg_i^t - \ybar^t}^2$, we find that
\begin{align}
    \norm{\yreg_j^t - \ybar^t} 
    &\leq \sum_{i=1}^n \norm{\yreg_i^t - \ybar^t} \leq \sqrt{n}\left( \sum_{i=1}^n \norm{\yreg_i^t - \ybar^t}^2 \right)^{1/2} \leq \frac{n^{3/2} D}{T(1-\beta)}
    \label{eq:conv_lemma_inequal}
\end{align}
where inequality (\ref{eq:conv_lemma_inequal}) follows from (\ref{eq:conv_lemma_final_bound}).
\end{proof}



\begin{lemma}
For any iteration $t\leq T$ for $T\in\Z_{++}$, the Euclidean distance between the means $\ybar^t$ and $\ybar^{t-1}$ of the local variables $\yreg_i^t$ and $\yreg_i^{t-1}$ respectively for $i\in\agents$ at consecutive iterations $t$ and $t-1$ can be bounded by
\begin{align}
    \norm{\ybar^{t} - \ybar^{t-1}}_2 \leq \frac{D}{T}.
    \label{eq:consensus_lemma_statement}
\end{align}

\end{lemma}\label{lemma:consensus}

\begin{proof}
Averaging both sides of the update rule for $\yreg_i^t$ of Algorithm \ref{algo:cdcg} across the set of agents $i\in\mathcal{N}$ yields the following expression for $\ybar^t$:
\begin{align}
    \ybar^t &= \frac{1}{n} \sum_{i=1}^n \sum_{j\in\neighbors} w_{ij} \yreg_j^{t-1} + \frac{1}{T} \sum_{i=1}^n \vreg_i^{t}.
    \label{eq:consensus_lemma_1}
\end{align}
Since $w_{ij}=0$ if $j\not\in\mathcal{N}_i\cup\{i\}$, we can rewrite the RHS of (\ref{eq:consensus_lemma_1}) in the following way:
\begin{align}
    \ybar^t &= 
    \frac{1}{n} \sum_{i=1}^n \sum_{j=1}^n w_{ij}\yreg_j^{t-1} + \frac{1}{T} \sum_{i=1}^n \vreg_i^{t} \notag \\
    &= \frac{1}{n} \sum_{j=1}^n \yreg_j^{t-1} \sum_{i=1}^n w_{ij} +  \frac{1}{T} \sum_{i=1}^n \vreg_i^{t} \notag \\
    &= \frac{1}{n} \sum_{j=1}^n \yreg_j^{t-1} + \frac{1}{T} \sum_{i=1}^n \vreg_i^{t}
    \label{eq:consensus_lemma_2}
\end{align}
where (\ref{eq:consensus_lemma_2}) follows since $\W^T\1 = \1$.  Rearranging (\ref{eq:consensus_lemma_2}), it follows that
\begin{align*}
    \norm{\ybar^t - \ybar^{t-1}} = \frac{1}{T} \norm{ \sum_{i=1}^n \vreg_i^{t}} \leq \frac{D}{T}
\end{align*}
Note that because the Euclidean distance between points of the polytope $P$ are assumed to be bounded, $\norm{\sum_{i=1}^n\vreg_i^{t}} \leq D$.  The expression in (\ref{eq:consensus_lemma_statement}) follows.
\end{proof}



\begin{cor}
Let $T\in\Z_{++}$.  Then the vector $\ybar = \frac{1}{n}\sum_{i=1}^n \yreg_i$ is in the constraint set $\setP$ $\forall t\leq T$.

\end{cor}\label{cor:feasibility}

\begin{proof}
In Lemma 1 we proved that $\yreg_i^t$ converges to $\ybar^t$.   We show that $\ybar^t \in\setP$ by induction.  Because we assign $\yreg_i^0 = \0$, it is clear that $\ybar^0\in\setP$.  Now as inductive hypothesis, we assume that $\ybar^{t-1}$ is in $\setP$.  Observe that we can write $\ybar^t = \ybar^{t-1} + (1/T)\sum_{i=1}^n \vreg_i^t$.  Thus by the inductive hypothesis and the fact that $\sum_{i=1}^n \vreg_i^t \in\setP$ $\forall t\leq T$, it follows that $\ybar^t$ is a convex combination of elements of $\setP$.  That is, we can write $\ybar^t = (1/T)\sum_{k=1}^t \sum_{i=1}^n \vreg_i^k + (1-t/T) \0$.  Therefore $\ybar^t\in\setP$, and so $\yreg_i^t$ converges to a point in $\setP$.
\end{proof}



\begin{lemma}
Let $F$ be the multilinear extension of a monotone submodular function $f:2^{\Y}\mapsto \R$ where $\Y$ is a discrete ground set.  Then
\begin{align}
    \max_{\vreg \in P_i \cap C_i} \langle\nabla F( \yreg_i), \vreg\rangle=\max_{\mathbf{x} \in P_i} \langle[\nabla F( \yreg_i)]_{c_i}, \mathbf{x}\rangle
    \label{eq:tech_lemma_argmax_statement}
\end{align}
where $[\nabla F(\bar \yreg_i)]_{c_i}$ denotes the projection of $\nabla F( \yreg_i)$ onto the set $C_i$.

\end{lemma}\label{lemma:tech_lemma_argmax}

\begin{proof}
Consider the definitions of $\setP_i$ and $\setC_i$ in (\ref{eq:def_p_i}) and (\ref{eq:def_c_i}) respectively.  Maximizing $\langle\nabla F( \yreg_i), \vreg\rangle$ over $\vreg\in\setP_i \cap \setC_i$ results in the same value as maximizing the inner product of the projection of $\nabla F(\yreg_i^{t-1})$ onto the set $\setC_i$ over $\xreg\in\setP_i$. 
\end{proof}



\begin{lemma}
Let $F$ be the multilinear extension of a monotone submodular function $f:2^{\Y}\mapsto \R$ where $\Y$ is a discrete ground set.  Then
\begin{align}
    \norm{\nabla F(\ybar^t)- \sum_{i=1}^{n}\left[\nabla F(\yreg_i^t) \right]_{C_i}}\leq \frac{n^{3/2}DL}{T(1-\beta)}
\end{align}

\end{lemma}\label{lemma:tech_lemma_bound}

\begin{proof}
Observe that
\begin{align}
  \norm{\nabla F(\ybar^t) - \sum_{i=1}^{n}\left[\nabla F(\yreg_i^t)\right]_{C_i}} &\leq \norm{\sum_{i=1}^{n} \left( \left[\nabla F(\ybar^t) \right]_{C_i}- \left[\nabla F(\yreg_i^t) \right]_{C_i} \right)} \notag \\
  &\leq \sum_{i=1}^n \norm{\left[ \nabla F(\ybar^t) \right]_{C_i} - \left[ \nabla F(\yreg_i^t) \right]_{C_i} } 
  \label{eq:tech_lem_bound_triangle_inequal} \\
  &\leq \sum_{i=1}^{n}\norm{\nabla F(\ybar^t) - \nabla F(\yreg_i^t)} 
  \label{eq:tech_lem_bound_proj} \\
  &\leq \frac{n^{3/2}DL}{T(1-\beta)}
  \label{eq:tech_lem_bound_result}
\end{align}
where (\ref{eq:tech_lem_bound_triangle_inequal}) follows from the triangle inequality, (\ref{eq:tech_lem_bound_proj}) follows by the definition of the set $\setC_i$, and (\ref{eq:tech_lem_bound_result}) follows from the assumption that $\nabla F$ is $L$-Lipschitz continuous and from Lemma \ref{lemma:convergence}.
\end{proof}

%% file: appendix/main-thm-proof.tex
\section*{Appendix C: Proof of Theorem \ref{thm:main_result}}


This Appendix establishes the main result of this paper, which is restated here for convenience.

\begin{theorem}
Consider the \texttt{CDCG} algorithm described in Algorithm \ref{algo:cdcg}.  Let $\yreg^*$ denote the global maximizer of the optimization problem defined in Problem \ref{prob:cont}, and assume that a positive integer $T$ and a doubly-stochastic weight matrix $\W$ are given.  Then provided that the assumptions outlined in Appendix A hold, for all nodes $i\in\mathcal{N}$, the local variables $\yreg_i^T$ obtained after $T$ iterations satisfy
\begin{align}
    F(\yreg_i^T) &\geq \left( 1 - \frac{1}{e} \right) F(\yreg^*) - \left[\frac{LD^2}{2T} + \frac{LD^2(n^2 + n^{5/2}) + n^{5/2}DG }{T(1-\beta) }\right] 
\end{align}
where $D$, $G$, $L$, and $\beta$ are problem-dependent constants that are formally defined in Appendices A and B.
\end{theorem}

\begin{proof}
Due to the assumption that $\nabla F$ is $L$-Lipschitz,
\begin{align}
    &F \left(\ybar^{t+1} \right) - F \left(\ybar^t \right) \notag \\ 
    &\quad\geq \left\langle \nabla F \left(\ybar^{t} \right), \ybar^{t+1} - \ybar^{t} \right\rangle \notag  - \frac{L}{2} \norm{\ybar^{t+1} - \ybar^t}^2 \notag \\
    &\quad \geq \left\langle \nabla F \left(\ybar^{t} \right), \ybar^{t+1} - \ybar^{t} \right\rangle - \frac{LD^2}{2T^2}
    \label{eq:main_thm_lipschitz}
\end{align}
where (\ref{eq:main_thm_lipschitz}) follows from Lemma \ref{lemma:consensus}.  Now consider that the inner-product term on the RHS of (\ref{eq:main_thm_lipschitz}) can be written in the following way:
\begin{align}
    \left\langle \nabla F \left(\ybar^{t} \right), \ybar^{t+1} - \ybar^{t} \right\rangle &=  \left\langle \nabla F \left(\ybar^{t}\right), \frac{1}{T}\sum_{i=1}^n \vreg_i^{t+1} \right\rangle \notag \\
    &= \frac{1}{T} \sum_{i=1}^n \Big[ \left\langle \nabla F(\ybar^t) - \nabla F(\yreg_i^t), \vreg_i^{t+1} \right\rangle  +  \left\langle \nabla F(\yreg_i^t), \vreg_i^{t+1} \right\rangle \Big].
    \label{eq:main_thm_inner_prod_term}
\end{align}
Here (\ref{eq:main_thm_inner_prod_term}) follows from the linearity of inner products and then from adding and subtracting $\nabla F(\yreg_i^t)$.  Our immediate goal is to bound (\ref{eq:main_thm_inner_prod_term}) from below.  To do so, consider that by the Cauchy-Schwartz inequality,
\begin{align}
    &\left\langle \nabla F \left(\ybar^{t} \right)-\nabla F(\yreg_i^t), \vreg_i^{t+1} \right\rangle \notag \\
    &\qquad\qquad\leq \norm{\nabla F(\ybar^t) - \nabla F(\yreg_i^t)} \cdot \norm{\vreg_i^{t+1}}
    \notag \\
    &\qquad\qquad\leq LD \norm{\ybar^t - \yreg_i^t}
    \label{eq:main_thm_lipschitz_bound} \\
    &\qquad\qquad\leq \frac{n^{3/2} LD^2}{T(1-\beta)}
    \label{eq:main_thm_conv_bound}
\end{align}
where (\ref{eq:main_thm_lipschitz_bound}) is due to the assumption that $\nabla F$ is $L$-Lipschitz continuous and (\ref{eq:main_thm_conv_bound}) follows from Lemma \ref{lemma:convergence}.  Next, because $\vreg_i^{t+1}$ is defined as the argmax between $\nabla F(\yreg_i^t)$ and vectors $\xreg\in \setP_i \cap \setC_i$ in the Step 3 of Algorithm \ref{algo:cdcg} and by Lemma \ref{lemma:tech_lemma_argmax} we have
\begin{align}
    \langle \nabla F(\yreg_i^t), \vreg_i^{t+1} \rangle 
    &\geq \langle [\nabla F(\yreg_i^t)]_{C_i}, \yreg^* \rangle 
    \label{eq:main_thm_proj_bound}.
\end{align}
By Lemma \ref{lemma:tech_lemma_bound}, if we let $\epsilon = \frac{n^{3/2}DL}{T(1-\beta)}$, we can conclude that
\begin{align}
    -\epsilon \1 + \nabla F(\ybar^t) \leq \sum_{i=1}^n \left[ \nabla F(\yreg_i^t) \right]_{\setC_i} \leq \nabla F(\ybar^t) + \epsilon\1.
    \label{eq:main_thm_eps_bound}
\end{align}
By construction, $\yreg^*\succeq 0$ since $\yreg^*\in\setP$.  Then we can infer from (\ref{eq:main_thm_eps_bound}) that
\begin{align}
    \left\langle \sum_{i=1}^n \left[ \nabla F(\yreg_i^t) \right]_{\setC_i} , \yreg^* \right\rangle \geq \left\langle -\epsilon \1, \yreg^* \right\rangle + \left\langle \nabla F(\ybar^t), \yreg^* \right\rangle.
    \label{eq:main_thm_one_way_eps_bound}
\end{align}
Our goal is to bound (\ref{eq:main_thm_one_way_eps_bound}).  To do this, consider that $\norm{\yreg^*} \leq D$ by (\ref{assume:unif-bounded}) and $\langle \1, \yreg^* \rangle = \norm{\yreg^*}_1$ since $\yreg^*\succeq 0$.  Since $\norm{\yreg^*}_1 \leq \sqrt{n}\norm{\yreg^*}_2$, we have $\langle \epsilon\1, \yreg^*\rangle \leq D\epsilon\sqrt{n}$.  Thus by replacing $\epsilon = \frac{n^{3/2}DL}{T(1-\beta)}$, we conclude that
\begin{align}
    \left\langle\sum_{i=1}^{n} \left[\nabla F(\yreg_i^t) \right]_{\setC_i}, \yreg^* \right\rangle 
    &\geq \langle \nabla F(\ybar^t), \yreg^* \rangle - \frac{n^2LD^2}{T(1-\beta)}\notag\\
    &\geq  F(\yreg^*)- F(\ybar^t) - \frac{n^2LD^2}{T(1-\beta)}
    \label{eq:main_thm_monotone_bound}.
\end{align}
Altogether, we have shown via (\ref{eq:main_thm_conv_bound}), (\ref{eq:main_thm_proj_bound}), and (\ref{eq:main_thm_monotone_bound}) that (\ref{eq:main_thm_inner_prod_term}) can be bounded by
\begin{align}
    &\langle \nabla F(\ybar^t), \bar{\yreg}^{t+1} - \bar{\yreg}^t \rangle \geq \frac{1}{T}  \left[ F(\yreg^*) - F(\ybar^t) - \frac{LD^2(n^2 + n^{5/2})}{T(1-\beta) } \right] 
    \label{eq:main_thm_first_term_bound}.
\end{align}
Furthermore, (\ref{eq:main_thm_first_term_bound}) and (\ref{eq:main_thm_lipschitz}) imply that
\begin{align}
    F \left(\ybar^{t+1} \right) - F \left(\ybar^t \right) &\geq \frac{1}{T} \left[ F(\yreg^*) - F(\ybar^t) \right] - \frac{LD^2(n^2 + n^{5/2})}{T^2(1-\beta) } - \frac{LD^2}{2T^2} \label{eq:main_thm_pre_rearrange}
\end{align}
Rearranging (\ref{eq:main_thm_pre_rearrange}), we obtain
\begin{align}
    F(\yreg^*) - F(\ybar^{t+1}) 
    &\leq \left( 1 - \frac{1}{T} \right) \left[ F(\yreg^*) - F(\ybar^t) \right] + \frac{LD^2(n^2 + n^{5/2})}{T^2(1-\beta) } + \frac{LD^2}{2T^2}
    \label{eq:main_thm_pre_reduce}.
\end{align}
By applying the inequality in (\ref{eq:main_thm_pre_reduce}) for $t=0, 1, \dots, T-1$, we find
\begin{align}
    &F(\yreg^*) - F(\ybar^T) \notag \\
    &\quad \leq \left( 1 - \frac{1}{T}\right)^T \left[ F(\yreg^*) - F(\ybar^0) \right] + \sum_{i=0}^{T-1} \left(1-\frac{1}{T}\right)^i \left[ \frac{LD^2(n^2 + n^{5/2})}{T^2(1-\beta) } + \frac{LD^2}{2T^2} \right] \notag  \\
    &\quad = \left(1-\frac{1}{T}\right)^T \left[ F(\yreg^*) - F(\ybar^0) \right] + \left(T - T\left(1 - \frac{1}{T}\right)^T\right)\left[\frac{LD^2(n^2 + n^{5/2})}{T^2(1-\beta) } + \frac{LD^2}{2T^2}\right] \notag \\
    &\quad\leq \frac{1}{e}\left[F(\yreg^*) - F(\ybar^0)\right]  + \left(1 - \frac{1}{e}\right)\left[\frac{LD^2(n^2 + n^{5/2})}{T(1-\beta) } + \frac{LD^2}{2T}\right] \notag \\
    &\quad\leq \frac{1}{e}\left[F(\yreg^*) - F(\ybar^0)\right]  + \left[\frac{LD^2(n^2 + n^{5/2})}{T(1-\beta) } + \frac{LD^2}{2T}\right]
    \label{eq:main_thm_common_denom}
\end{align}
where to derive (\ref{eq:main_thm_common_denom}) we used $(1-1/T)^T \leq 1/e$.  Now recall that we set $\yreg_i^0 = \0$.  Then from equation \eqref{eq:cont_ext_obj_func}, we have $F(\0) \geq 0$ $\forall i\in\agents$.  Thus follows that
\begin{align}
    F(\ybar^T) &\geq \left( 1 - \frac{1}{e} \right) F(\yreg^*) - \left[ \frac{LD^2(n^2 + n^{5/2})}{T(1-\beta) } + \frac{LD^2}{2T}\right].
    \label{eq:main_thm_final_bound}
\end{align}
Now by the assumption made in (\ref{assume:bounded-grad}), $F$ is $G$-Lipschitz continuous and therefore
\begin{align}
    \left| F(\ybar^T) - F(\yreg_i^T) \right| 
    &\leq G\norm{\ybar^T - \yreg_i^T} \leq \frac{n^{3/2}DG}{T(1-\beta)}
    \label{eq:main_thm_mean_to_pt_bound}
\end{align}
where (\ref{eq:main_thm_mean_to_pt_bound}) follows from Lemma \ref{lemma:convergence}.  Thus by combining the results in (\ref{eq:main_thm_final_bound}) and (\ref{eq:main_thm_mean_to_pt_bound}) we find that $\forall i\in\agents$, 
\begin{align*}
    F(\yreg_i^T) &\geq \left( 1 - \frac{1}{e} \right) F(\yreg^*) - \left[ \frac{LD^2}{2T} + \frac{LD^2(n^2 + n^{5/2}) + n^{5/2}DG }{T(1-\beta) }\right]
\end{align*}
and the claim in (\ref{thm:main_result_statement}) follows.
\end{proof}

%% file: l4dc2021-sample.bbl
\begin{thebibliography}{31}
\providecommand{\natexlab}[1]{#1}
\providecommand{\url}[1]{\texttt{#1}}
\expandafter\ifx\csname urlstyle\endcsname\relax
  \providecommand{\doi}[1]{doi: #1}\else
  \providecommand{\doi}{doi: \begingroup \urlstyle{rm}\Url}\fi

\bibitem[Adibi et~al.(2020)Adibi, Mokhtari, and Hassani]{adibi2020submodular}
Arman Adibi, Aryan Mokhtari, and Hamed Hassani.
\newblock Submodular meta-learning.
\newblock \emph{Advances in Neural Information Processing Systems}, 33, 2020.

\bibitem[Barbosa et~al.(2015)Barbosa, Ene, Nguyen, and Ward]{barbosa2015power}
Rafael Barbosa, Alina Ene, Huy Nguyen, and Justin Ward.
\newblock The power of randomization: Distributed submodular maximization on
  massive datasets.
\newblock In \emph{International Conference on Machine Learning}, pages
  1236--1244, 2015.

\bibitem[Buchbinder et~al.(2014)Buchbinder, Feldman, Naor, and
  Schwartz]{buchbinder2014submodular}
Niv Buchbinder, Moran Feldman, Joseph Naor, and Roy Schwartz.
\newblock Submodular maximization with cardinality constraints.
\newblock In \emph{Proceedings of the twenty-fifth annual ACM-SIAM symposium on
  Discrete algorithms}, pages 1433--1452. SIAM, 2014.

\bibitem[Buchbinder et~al.(2015)Buchbinder, Feldman, Seffi, and
  Schwartz]{buchbinder2015tight}
Niv Buchbinder, Moran Feldman, Joseph Seffi, and Roy Schwartz.
\newblock A tight linear time (1/2)-approximation for unconstrained submodular
  maximization.
\newblock \emph{SIAM Journal on Computing}, 44\penalty0 (5):\penalty0
  1384--1402, 2015.

\bibitem[Calinescu et~al.(2011)Calinescu, Chekuri, Pal, and
  Vondr{\'a}k]{calinescu2011maximizing}
Gruia Calinescu, Chandra Chekuri, Martin Pal, and Jan Vondr{\'a}k.
\newblock Maximizing a monotone submodular function subject to a matroid
  constraint.
\newblock \emph{SIAM Journal on Computing}, 40\penalty0 (6):\penalty0
  1740--1766, 2011.

\bibitem[Chekuri et~al.(2014)Chekuri, Vondr{\'a}k, and
  Zenklusen]{chekuri2014submodular}
Chandra Chekuri, Jan Vondr{\'a}k, and Rico Zenklusen.
\newblock Submodular function maximization via the multilinear relaxation and
  contention resolution schemes.
\newblock \emph{SIAM Journal on Computing}, 43\penalty0 (6):\penalty0
  1831--1879, 2014.

\bibitem[Chen et~al.(2020)Chen, Zhang, Hassani, and Karbasi]{chen2020black}
Lin Chen, Mingrui Zhang, Hamed Hassani, and Amin Karbasi.
\newblock Black box submodular maximization: Discrete and continuous settings.
\newblock In \emph{International Conference on Artificial Intelligence and
  Statistics}, pages 1058--1070, 2020.

\bibitem[Clark et~al.(2015)Clark, Alomair, Bushnell, and
  Poovendran]{clark2015scalable}
Andrew Clark, Basel Alomair, Linda Bushnell, and Radha Poovendran.
\newblock Scalable and distributed submodular maximization with matroid
  constraints.
\newblock In \emph{2015 13th International Symposium on Modeling and
  Optimization in Mobile, Ad Hoc, and Wireless Networks (WiOpt)}, pages
  435--442. IEEE, 2015.

\bibitem[Corah and Michael(2017)]{corah2017efficient}
Micah Corah and Nathan Michael.
\newblock Efficient online multi-robot exploration via distributed sequential
  greedy assignment.
\newblock In \emph{Robotics: Science and Systems}, volume~13, 2017.

\bibitem[Corah and Michael(2018)]{corah2018distributed}
Micah Corah and Nathan Michael.
\newblock Distributed submodular maximization on partition matroids for
  planning on large sensor networks.
\newblock In \emph{2018 IEEE Conference on Decision and Control (CDC)}, pages
  6792--6799. IEEE, 2018.

\bibitem[Corah and Michael(2019)]{corah2019distributed}
Micah Corah and Nathan Michael.
\newblock Distributed matroid-constrained submodular maximization for
  multi-robot exploration: Theory and practice.
\newblock \emph{Autonomous Robots}, 43\penalty0 (2):\penalty0 485--501, 2019.

\bibitem[Djolonga et~al.(2016)Djolonga, Tschiatschek, and
  Krause]{djolonga2016variational}
Josip Djolonga, Sebastian Tschiatschek, and Andreas Krause.
\newblock Variational inference in mixed probabilistic submodular models.
\newblock In \emph{Advances in Neural Information Processing Systems}, pages
  1759--1767, 2016.

\bibitem[Du et~al.(2020)Du, Qian, Claudel, and Sun]{du2020jacobi}
Bin Du, Kun Qian, Christian Claudel, and Dengfeng Sun.
\newblock Jacobi-style iteration for distributed submodular maximization.
\newblock \emph{arXiv preprint arXiv:2010.14082}, 2020.

\bibitem[Gharesifard and Smith(2017)]{gharesifard2017distributed}
Bahman Gharesifard and Stephen~L Smith.
\newblock Distributed submodular maximization with limited information.
\newblock \emph{IEEE transactions on control of network systems}, 5\penalty0
  (4):\penalty0 1635--1645, 2017.

\bibitem[Golovin and Krause(2011)]{golovin2011adaptive}
Daniel Golovin and Andreas Krause.
\newblock Adaptive submodularity: Theory and applications in active learning
  and stochastic optimization.
\newblock \emph{Journal of Artificial Intelligence Research}, 42:\penalty0
  427--486, 2011.

\bibitem[Grimsman et~al.(2018)Grimsman, Ali, Hespanha, and
  Marden]{grimsman2018impact}
David Grimsman, Mohd~Shabbir Ali, Joao~P Hespanha, and Jason~R Marden.
\newblock The impact of information in greedy submodular maximization.
\newblock \emph{IEEE Transactions on Control of Network Systems}, 2018.

\bibitem[Hassani et~al.(2017)Hassani, Soltanolkotabi, and
  Karbasi]{hassani2017gradient}
Hamed Hassani, Mahdi Soltanolkotabi, and Amin Karbasi.
\newblock Gradient methods for submodular maximization.
\newblock In \emph{Advances in Neural Information Processing Systems}, pages
  5841--5851, 2017.

\bibitem[Hu et~al.(2007)Hu, Chen, Lou, and Li]{hu2007distributed}
Yusuo Hu, Hua Chen, Jian-guang Lou, and Jiang Li.
\newblock Distributed density estimation using non-parametric statistics.
\newblock In \emph{27th International Conference on Distributed Computing
  Systems (ICDCS'07)}, pages 28--28. IEEE, 2007.

\bibitem[Mirzasoleiman et~al.(2013)Mirzasoleiman, Karbasi, Sarkar, and
  Krause]{mirzasoleiman2013distributed}
Baharan Mirzasoleiman, Amin Karbasi, Rik Sarkar, and Andreas Krause.
\newblock Distributed submodular maximization: Identifying representative
  elements in massive data.
\newblock In \emph{Advances in Neural Information Processing Systems}, pages
  2049--2057, 2013.

\bibitem[Mirzasoleiman et~al.(2016)Mirzasoleiman, Karbasi, Sarkar, and
  Krause]{mirzasoleiman2016distributed}
Baharan Mirzasoleiman, Amin Karbasi, Rik Sarkar, and Andreas Krause.
\newblock Distributed submodular maximization.
\newblock \emph{The Journal of Machine Learning Research}, 17\penalty0
  (1):\penalty0 8330--8373, 2016.

\bibitem[Mokhtari et~al.(2018)Mokhtari, Hassani, and
  Karbasi]{mokhtari2018decentralized}
Aryan Mokhtari, Hamed Hassani, and Amin Karbasi.
\newblock Decentralized submodular maximization: Bridging discrete and
  continuous settings.
\newblock \emph{arXiv preprint arXiv:1802.03825}, 2018.

\bibitem[Mokhtari et~al.(2020)Mokhtari, Hassani, and
  Karbasi]{mokhtari2020stochastic}
Aryan Mokhtari, Hamed Hassani, and Amin Karbasi.
\newblock Stochastic conditional gradient methods: From convex minimization to
  submodular maximization.
\newblock \emph{Journal of Machine Learning Research}, 21\penalty0
  (105):\penalty0 1--49, 2020.

\bibitem[Nemhauser and Wolsey(1978)]{nemhauser1978best}
George~L Nemhauser and Laurence~A Wolsey.
\newblock Best algorithms for approximating the maximum of a submodular set
  function.
\newblock \emph{Mathematics of operations research}, 3\penalty0 (3):\penalty0
  177--188, 1978.

\bibitem[Nemhauser et~al.(1978)Nemhauser, Wolsey, and
  Fisher]{nemhauser1978analysis}
George~L Nemhauser, Laurence~A Wolsey, and Marshall~L Fisher.
\newblock An analysis of approximations for maximizing submodular set
  functions—i.
\newblock \emph{Mathematical programming}, 14\penalty0 (1):\penalty0 265--294,
  1978.

\bibitem[Schlotfeldt et~al.(2018)Schlotfeldt, Thakur, Atanasov, Kumar, and
  Pappas]{schlotfeldt2018anytime}
Brent Schlotfeldt, Dinesh Thakur, Nikolay Atanasov, Vijay Kumar, and George~J
  Pappas.
\newblock Anytime planning for decentralized multirobot active information
  gathering.
\newblock \emph{IEEE Robotics and Automation Letters}, 3\penalty0 (2):\penalty0
  1025--1032, 2018.

\bibitem[Singh et~al.(2009)Singh, Krause, Guestrin, and
  Kaiser]{singh2009efficient}
Amarjeet Singh, Andreas Krause, Carlos Guestrin, and William~J Kaiser.
\newblock Efficient informative sensing using multiple robots.
\newblock \emph{Journal of Artificial Intelligence Research}, 34:\penalty0
  707--755, 2009.

\bibitem[Wei et~al.(2013)Wei, Liu, Kirchhoff, and Bilmes]{wei2013using}
Kai Wei, Yuzong Liu, Katrin Kirchhoff, and Jeff Bilmes.
\newblock Using document summarization techniques for speech data subset
  selection.
\newblock In \emph{Proceedings of the 2013 Conference of the North American
  Chapter of the Association for Computational Linguistics: Human Language
  Technologies}, pages 721--726, 2013.

\bibitem[Wolsey(1982)]{wolsey1982analysis}
Laurence~A Wolsey.
\newblock An analysis of the greedy algorithm for the submodular set covering
  problem.
\newblock \emph{Combinatorica}, 2\penalty0 (4):\penalty0 385--393, 1982.

\bibitem[Xie et~al.(2019)Xie, Zhang, Shen, Mi, and Qian]{xie2019decentralized}
Jiahao Xie, Chao Zhang, Zebang Shen, Chao Mi, and Hui Qian.
\newblock Decentralized gradient tracking for continuous dr-submodular
  maximization.
\newblock In \emph{The 22nd International Conference on Artificial Intelligence
  and Statistics}, pages 2897--2906, 2019.

\bibitem[Zhong and Cassandras(2011)]{zhong2011distributed}
Minyi Zhong and Christos~G Cassandras.
\newblock Distributed coverage control and data collection with mobile sensor
  networks.
\newblock \emph{IEEE Transactions on Automatic Control}, 56\penalty0
  (10):\penalty0 2445--2455, 2011.

\bibitem[Zhou et~al.(2020)Zhou, Tzoumas, Pappas, and
  Tokekar]{zhou2020distributed}
Lifeng Zhou, Vasileios Tzoumas, George~J Pappas, and Pratap Tokekar.
\newblock Distributed attack-robust submodular maximization for multi-robot
  planning.
\newblock In \emph{2020 IEEE International Conference on Robotics and
  Automation (ICRA)}, pages 2479--2485. IEEE, 2020.

\end{thebibliography}
